\documentclass{amsart}

\newcommand{\ga}{\alpha}
\newcommand{\gb}{\beta}
\newcommand{\ggm}{\gamma}
\newcommand{\gd}{\delta}

\newcommand{\eps}{\varepsilon}

\newtheorem{lemma}{Lemma}
\renewcommand{\epsilon}{\varepsilon}
\renewcommand{\phi}{\varphi}
\newtheorem{remark}{Remark}
\newtheorem{theorem}{Theorem}
\newtheorem{prop}{Proposition}

\newcommand{\pa}{\partial}

\begin{document}
\parindent0 in
\parskip 1 em
\title{On the Geometry of Moduli Space of Polarized Calabi-Yau manifolds}
\author{Michael Douglas and Zhiqin Lu}
\thanks{The first author is partially supported by DOE grant DE-FG02-96ER40959, and 
the second
author is partially supported by  NSF Career award DMS-0347033 and the
Alfred P. Sloan Research Fellowship.}
\email[Michael Douglas, NHETC and Department of Physics and Astronomy, Rutgers University, Piscataway, NJ 08855-0849, USA and IHES, Le Bois-Marie, Bures-sur-Yvette, 91440, France]{mrd@physics.rutgers.edu}

\email[Zhiqin Lu, Department of Mathematics, UC Irvine, Irvine, CA 92697, USA]{zlu@math.uci.edu}

\date{December 14, 2005}

\maketitle
\section{Introduction}

Let $X$ be a compact K\"ahler manifold with zero first Chern class,
and let $L$ be an ample line bundle over $X$. The pair $(X,L)$ is
called a polarized Calabi-Yau manifold.  By Yau's proof of the Calabi
conjecture, we know such a manifold carries a unique Ricci flat metric
compatable with the polarization (cf. ~\cite{y1}).  Thus, the moduli
space of such Ricci flat K\"ahler metrics is
the moduli space of complex structures of $(X,L)$.

By a theorem of Mumford, a Calabi-Yau moduli space (or any coarse moduli
space of polarized K\"ahler manifolds) is a complex variety. By Tian~\cite{t2}, locally, up to a finite cover, the moduli space is smooth.  Now,
in Riemannian geometry, there is a natural metric on any moduli space
of metrics, the Weil-Petersson metric, obtained by restriction from
the metric on the space of metrics.  Quite a lot is known about the
local structure of the WP metric on Calabi-Yau moduli space.  But much
less is known about its global properties.

In this short paper, we study the integrals of the curvature
invariants of the Weil-Petersson metric on a Calabi-Yau moduli space.
In Theorem~\ref{main1}, we prove that these quantities are all finite.
In Theorem~\ref{main2} and in work to appear~\cite{DLN}, we prove that
they are rational numbers.  Now if the moduli space had been compact,
then this would be expected by the theorem of Gauss-Bonnet-Chern.  But
Calabi-Yau moduli spaces are not compact, making this result
nontrivial.

Besides its mathematical interest, the geometry of Calabi-Yau moduli
space is very interesting in string theory, and there are various
physics arguments~\cite{HorneMoore,Ad,strings,DL, vafa} suggesting the
finiteness of the volume and integrability of the curvature invariants
of the Weil-Petersson metric.  Recently Eguchi and Tachikawa have
shown the finiteness of the curvature invariant of \cite{Ad} near a
large subclass of singularities of moduli space \cite{Eguchi:2005eh},
calling upon mathematical methods similar to those we use below.

Mathematically, this paper is a continuation of the previous works in
~\cite{lu3, lu5,lu12, ls-1, ls-2, fl1,fly1}, on the local and global
geometry of the moduli space and the BCOV torsion of Calabi-Yau
moduli.

Before finishing this section, we write out explicitly the Calabi-Yau
moduli of the most famous Calabi-Yau threefold: the quintic
hypersurface in $CP^4$.  Let this be
\[
X=\{Z\mid Z_0^5+\cdots+Z_4^5+5\lambda Z_0\cdots Z_4=0\}\subset CP^4.
\]
It is a smooth hypersurface if $\lambda$ is not any of the fifth unit
roots. To construct the moduli space, we define
\[
V=\{f\mid f \text{ is a homogeneous quintic polynomial of $Z_0,\cdots,Z_4$}\}.
\]
one can verify that $\dim V=126$. Thus for any $t\in P(V)=CP^{125}$,
$t$ is represented by a hypersurface. However, if two hypersurfaces
differ by an element in $Aut(CP^4)$, then they are considered the
same. Let $D$ be the divisor in $CP^{125}$ characterizing the singular
hypersurfaces in $CP^4$. Then the moduli space of $X$ is
\[
\mathcal M=CP^{125}\backslash D/Aut(CP^4).
\]
The dimension of the moduli space is $101$. But other than the
dimension, we still know very little about this variety.

The organization of the paper is as follows: in Section 2, we give
some physics background of our problems; in Section 3, we define the
Weil-Petersson metric; in Section 4, we present the main results of
this paper; then we introduce the Hodge metrics in Section 5; in the
last section, we prove the main results of this paper.

{\bf Acknowledgement.} The second author presented the main results of
this paper in the 2005 RIMS Symposium on {\it Analytic Geometry of the
Bergman kernel and Related Topics}. He thanks the orgainzers,
especially Professor Ohsawa, for the hospitality during his visit of
RIMS.

\section{Physics background}

In the original compactifications of heterotic string theory
\cite{GSW}, as well as in many later constructions, the universe is a
direct product of a $4d$ space-time and a tiny, compact Ricci flat six
manifold $M$.  Arguments from supersymmetry, as well as the fact that
we know no other examples, suggest that $M$ is a Calabi-Yau manifold.
While we do not know which $M$ to choose, we do know how to go from
geometric properties of $M$, together with certain auxiliary data, to
statements about observable physics.  Then, if a particular choice of
$M$ and the auxiliary data implies statements which are in conflict
with observation, we know this choice is incorrect.  At present it is
an open problem to show that any specific choice or ``vacuum'' is
consistent with current observations.  Given that such choices exist,
we would like to go on to show that the number of vacua is finite, and
estimate their number.

Suppose we assume a particular Calabi-Yau $M$; then the fact that
Ricci flat metrics on $M$ come in moduli spaces leads to the existence
of approximate solutions, in which the moduli of $M$ are slowly
varying in four-dimensional space-time.  These lead almost inevitably
to corrections to Newton's (and Einstein's) laws of gravity which
contradict observation, and thus we must somehow modify the construction
by postulating additional background fields to remove these moduli.

One way to do this is flux compactification, described in \cite{encyc}
and references there.  This construction picks out special points in
moduli space, the flux vacua, and thus part of counting vacua is to count
these points.

In their studies of type IIb flux compactification, the first author
and his collaborators (cf. ~\cite{Ad,Dou, DSZ-1,DSZ-2,DSZ-3}) derived
an asymptotic formula for the number of flux vacua,
\cite[Theorem 1.8]{DSZ-3}.  It is a product of a coefficient determined
by topological data of $M$, with an integral of a curvature invariant
derived from the Weil-Petersson metric over the moduli space.  Thus, the
finiteness of such integrals implies the finiteness of the number of
flux vacua (up to certain caveats explained in \cite{DSZ-3}), which was
a primary motivation for us to write this paper.

We are also very interested to the duality between of the special
K\"ahler manifolds and the Calabi-Yau moduli. In the proof of the
finiteness of the volume of Calabi-Yau moduli~\cite{ls-1} and in the
proof of the incompleteness of special K\"ahler manifold~\cite{lu8},
we use the generalized maximal principal. It would be interesting to
answer the following questions: Are the volume of projective special
K\"ahler manifolds finite? Are the Calabi-Yau moduli always incomplete
with respect to the Weil-Petersson metric? We hope that not only one
can answer these questions but also we can find relations of these two
problems.

\section{Weil-Petersson geometry}

Let ${\mathcal M}$ be the moduli space of a polarized Calabi-Yau
manifold of dimension $n$. Let $0\subset F^n\subset
F^{n-1}\subset\cdots \subset F^1\subset F^0= H$ be the Hodge bundles over
$\mathcal M$. Since each point of $\mathcal M$ is represented by a
Calabi-Yau manifold, the rank of $F^n$ is $1$. A natural Hermitian
metric on $F^n$ is given by the second Hodge-Riemann relation:
\[
C\int_{X_t}\Omega\wedge\bar\Omega>0,
\]
where $C$ is suitable constant and $\Omega$ is a nonzero $(n,0)$ form
of $X_t, t\in{\mathcal M}$. By a theorem of Tian~\cite{t2}, we know
that the curvature of the above Hermitian metric is positive, and the
Weil-Petersson metric is equal to the curvature of $F^n$. Thus we can
{\sl define} the Weil-Petersson metric whose K\"ahler form is the
curvature form of the line bundle. Let the K\"ahler form of the
Weil-Petersson metric be $\omega_{WP}$, then we have
\[
\omega_{WP}=-\sqrt{-1}\pa\bar\pa\log \int_{X_t}\Omega\wedge\bar\Omega,
\]
where $\Omega$ is a local holomorphic section of the bundle $F^n$.

The Weil-Petersson geometry is composed of the moduli space $\mathcal
M$, the Hodge bundles $F^k, k=0,\cdots,n$, and the Weil-Petersson
metric $\omega_{WP}$.  In order to understand the geometry of the
moduli space, we need to study the
curvature and the asymptotic behavior of the Weil-Petersson metric. Let $(\,,\,)$ be
the quadratic form on $H$ defined by the cup product. The quadratic
form is nondegenerate but not positive definite. Let
\[
F^k=H^{n,0}\oplus\cdots \oplus H^{k,n-k}
\]
be the orthogonal splitting with respect to the quadratic form
$(\,,\,)$.  Let $(\frac{\pa}{\pa t_1},\cdots, \frac{\pa}{\pa t_m})$ be
a local holomorphic frame near a smooth point $x$ of $\mathcal
M$. Define $\nabla_i\Omega$ to be the $H^{n-1,1}$ part of
$\pa_i\Omega=\frac{\pa\Omega}{\pa t_i}$ and $\nabla_j\nabla_i\Omega$
to be the $H^{n-2,2}$ part of $\pa_j\pa_i\Omega$ or
$\pa_j\nabla_i\Omega$. Then we have the following result:

\begin{theorem} \label{thm1}The curvature tensor $R_{\alpha\bar\beta\gamma\bar\delta}$ of the Weil-Petersson metric is
\begin{equation}\label{str1}
R_{\alpha\bar\beta\gamma\bar\delta}
=g_{\ga\bar\gb}g_{\ggm\bar\gd}+g_{\ga\bar\gd}g_{\ggm\bar\gb}-
\frac{(\nabla_\ga \nabla_\ggm\Omega,\overline{\nabla}_\gb \overline{\nabla}_\gd\Omega)}{(\Omega,\bar\Omega)}.
\end{equation}
\end{theorem}

\qed

If $n=3$, then we have
\[
\frac{(\nabla_\ga \nabla_\ggm\Omega,\overline{\nabla}_\gb \overline{\nabla}_\gd\Omega)}{(\Omega,\bar\Omega)}=F_{\alpha\gamma m}\overline{F_{\beta\delta n}}g^{m\bar n}/(\Omega,\bar\Omega)^2,
\]
where $\{F_{\alpha\beta\ggm}\}$ is the Yukawa coupling, which is a holomorphic section of the bundle ${\rm Sym}^{\otimes 2}\, F^n\otimes{\rm Sym}^{\otimes 3} T^*{\mathcal M}$, locally   defined as
\[
F_{\alpha\beta\gamma}=(\Omega,\pa_\alpha\pa_\beta\pa_\gamma\Omega).
\]
So ~\eqref{str1} can be written as 
\begin{equation}\label{str}
R_{\alpha\bar\beta\gamma\bar\delta}
=g_{\ga\bar\gb}g_{\ggm\bar\gd}+g_{\ga\bar\gd}g_{\ggm\bar\gb}-
F_{\alpha\gamma m}\overline{F_{\beta\delta n}}g^{m\bar n}/(\Omega,\bar\Omega)^2.
\end{equation}

\begin{remark} 
Formula~\eqref{str} was first given by Strominger~\cite{s}. In the
general case, a Hodge theoretic proof of Theorem~\ref{thm1} was given
by Wang~\cite{wang1}. Schumacher's paper ~\cite{sc1} will lead another
proof using the method of Siu~\cite{siu2}.
\end{remark}

We know very little of the global behavior of the moduli space
$\mathcal M$, except the following result of Viehweg~\cite{v2}.

\begin{theorem}[Viehweg] $\mathcal M$ is quasi-projective.
\end{theorem}

\begin{remark}\label{rem1}
By the above theorem, after normalization and desingularization, there
is a compact manifold $\overline{\mathcal M}$ such that
$\overline{\mathcal M}\backslash{\mathcal M}$ is a divisor of normal
crossings. In fact, since in general $\mathcal M$ is a complex
variety, we can redefine $\mathcal M$ to be the regular part of
$\mathcal M$. On such a setting, up to a finite cover, both $\mathcal
M$ and $\overline{\mathcal M}$ are manifolds.
\end{remark}

For the extension of the Hodge bundles across the divisor at infinity,
we have the following theorem of Schmid~\cite{schmid} or
Steenbrink~\cite{steenbrink}:

\begin{theorem} 
Let $\mathfrak X\rightarrow\Delta^r\times(\Delta^*)^s$ be a family of
 polarized Calabi-Yau manifolds, where $\Delta$ and $\Delta^*$ are the
 unit disk and the punctured unit disk, respectively. Suppose that all
 the monodromy operators are unipotent. Then there is a natural
 extension of the Hodge bundles to $\Delta^{r+s}$.
\end{theorem}

By the following result on the Weil-Petersson metric, the extension of
the bundles $F^n, F^{n-1}$ will give us information of the limiting
behaviors of the Weil-Petersson metric at infinity.

We need the following result of Tian~\cite{t2} in the rest of this paper:

\begin{prop}\label{pp1}
Let $(g_{\alpha\bar\beta})$ be the metric matrix of the Weil-Petersson
metric under the frame $(\frac{\pa}{\pa t_1},\cdots, \frac{\pa}{\pa
t_m})$. Then we have
\[
g_{\alpha\bar\beta}=-\frac{(\nabla_\alpha\Omega,\overline{\nabla_\beta\Omega})}{(\Omega,\overline\Omega)}.
\]
\end{prop}

The proof is a straightforward computation and is omitted.

\qed

\section{The main results}
There are several previous results related to the main results of this
paper.  It was proved in~\cite[Theorem 5.2]{ls-1} that the volume with
respect to the Weil-Petersson metric and the Hodge metric is
finite. In~\cite{ls-2, todo-1}, the rationality of the volume with respect to
the Weil-Peterssoin metric was proved. Furthermore, in~\cite{ls-2}, it was proved that the
integration of the $n$-th power of the Ricci curvature of the
Weil-Petersson metric is a rational number. The main results of this
paper are Theorem~\ref{main1} and Theorem~\ref{main2}. The most
general forms of these results will appear in our upcoming
paper~\cite{DLN}.

\begin{theorem} \label{main1}
Let $R_{WP}$ be the curvature tensor of $\omega_{WP}$. Let
$R=R_{WP}\otimes 1+1\otimes\omega_{WP}$. Let $f$ be any invariant
polynomial of $R$. Then we have
\[
\int_{\mathcal M} f(R)<+\infty.
\]
\end{theorem}

The above theorem is equivalent to the following: let $f_1,\cdots,f_s$
be invariant polynomials of $R_{WP}$ of degree $k_1,\cdots,k_s$,
respectively. Then
\[
\int_{\mathcal M}\sum_i f_i(R_{WP})\wedge\omega_{WP}^{m-k_i}<+\infty,
\]
where $m$ is the complex dimension of $\mathcal M$. The result is a generalization of Theorem 5.2 in~\cite{ls-1}.

\begin{theorem}\label{main2}
We assume that $\dim\,\mathcal M=2$.  Let $R$ be the curvature
operator of the Hodge bundle, and let $f$ be an invariant polynomial
with rational coefficients. Then we have
\[
\int_{\mathcal M} f(R)\in\mathbb Q.
\]
\end{theorem}

If $\dim\mathcal M=1$, or the rank of $F^k$ is one, then the
corresponding result follows from~\cite{ls-2}, because the only Chern
class will be the first Chern class. The case that $\mathcal M$ is of
arbitrary dimension is treated in~\cite{DLN}.

\section{The Hodge metric}

The curvature properties of the Weil-Petersson metric are not
good. For example, even in the case when the moduli space is of
dimension $1$, from~\cite[page 65]{cogp}, we know that the sign of the
Gauss curvature is not fixed. In~\cite{lu3}, the second author
introduced another natural metric, called the Hodge metric, on
$\mathcal M$.  We shall see that the Hodge metric is the bridge
between the curvature invariants and finiteness.

The following definition of the Hodge metric is from~\cite[section
6]{ls-1}, which is slightly different from that in~\cite{lu3}.

Let $\mathcal M$ be the moduli space of any polarized compact Kahler
manifold (not necessarily Calabi-Yau). Let $x\in\mathcal M$ be a
smooth point of $\mathcal M$. Assume that the period map $p:\mathcal
M\rightarrow D$ is an immersion near $x$.

Let $0\subset F^n\subset\cdots \subset F^1\subset F^0=H$ be the Hodge
bundles and let
\[
T_t{\mathcal M}\rightarrow H^1(X_t,\Theta_t)
\]
be the Kodaira-Spencer isomorphism. Let $\xi\in T_t{\mathcal M}$. Then $\xi$ defines a map
\[
\xi: H^1(X_t,\Theta_t)\times H^{p,q}\rightarrow H^{p-1,q+1}.
\]
Let $||\xi||_{p,q}$  be the operator norm with respect to the metric on Hodge bundles $H^{p,q}$ and $H^{p-1,q+1}$. Then we define
\[
||\xi||^2=\sum_{p+q=n}||\xi||^2_{p,q}.
\]
From the above definition, we get a Hermitian metric on the smooth
part of the moduli space $\mathcal M$. Let $\omega_H$ be the K\"ahler
form of the metric. Then the properties of the Hodge metric can be
summarized as follows~\cite{lu3}:

\begin{theorem}
Using the above notations, we have
\begin{enumerate}
\item The Hodge metric is a K\"ahler metric;
\item The bisectional curvartures of the Hodge metric are nonpositive;
\item The holomorphic sectional curvatures of the Hodge metric are  bounded from above by a negative constant;
\item The Ricci curvature of the Hodge metric  is bounded above by a negative constant.
\end{enumerate}
\end{theorem}

\begin{remark} If $p$ is not an immersion, then using the same definition, we get a 
semi-positive pseudo metric and the form $\omega_H$ is also well
defined. In fact, up to a constant, the Hodge metric is the pull-back
of the invariant Hermitian metric of $D$, the classifying space.  Such
a metric is K\"ahler in the sense that $d\omega_H=0$.  The pseudo
metric is called a generalized Hodge metric in~\cite{fl1}.
\end{remark}

The Hodge metric and the Weil-Petesson metric
on the moduli space of polarized Calabi-Yau manifolds
are closely related.
Let the dimension of the moduli space be $m$. Then we have the following
\begin{theorem} Using the above notations, we have
\begin{enumerate}
\item By Proposition~\ref{pp1}, we have
\[
2\omega_{WP}\leq\omega_{H}.
\]
\item If $\mathcal M$ is the moduli space of algebraic K3 surfaces, then
\[
2\omega_{WP}=\omega_H.
\]
\item If $\mathcal M$ is the moduli space of a Calabi-Yau threefold, then
we have~\cite{lu5}
\[
\omega_H=(m+3)\omega_{WP}+{\rm Ric}(\omega_{WP}).
\]
\item If $\mathcal M$ is the moduli space of a Calabi-Yau fourfold, then
we have~\cite{ls-1}
\[
\omega_H=2(m+2)\omega_{WP}+2{\rm Ric}(\omega_{WP}).
\]
\item In general, the generalized Hodge metrics and the Weil-Petersson
metric are related by the so-called BCOV torsion
(cf. ~\cite{bcov2,bcov1}) in~\cite{fl1}:
\[
\sum_{i=1}^n(-1)^i\omega_{H^i}-\frac{\sqrt{-1}}{2\pi}\pa\bar\pa\log T=\frac{\chi}{12}\omega_{WP},
\]
where $T$ is the BCOV torsion, $\omega_{H^i}$ is the generalized Hodge
metric with respect to the variation of Hodge structures of weight
$i$, and $\chi$ is the Euler characteristic number of a generic fiber.
\end{enumerate}
\end{theorem}

In order to study the asymptotic behavior of the Hodge metrics, we
quote the following Schwarz lemma of Yau~\cite{y3}:

\begin{theorem}\label{yau} Let $M,N$ be K\"ahler
manifolds. Suppose that $M$ is complete and the Ricci curvature of $M$
is bounded from below. Suppose that the bisectional curvatures of $N$
are nonpositive and the holomorphic sectional curvatures are bounded
from above by a negative constant. Then there is a constant $C$,
depending only on the dimensions of the two manifolds and the above
curvature bounds, such that
\[
f^*(\omega_N)\leq C\omega_M,
\]
where $\omega_M,\omega_N$ are the K\"ahler forms of manifold $M$ and
$N$, respectively.
\end{theorem}

By Remark~\ref{rem1}, the regular part of the complex variety
$\mathcal M$ is quasi-projective.  We construct a K\"ahler metric
$\omega_P$ on a quasi-projective manifold, following
Jost-Yau~\cite{jy1}. Let $U=(\Delta^*)^r\times\Delta^s$. We define a
K\"ahler metric on $U$ by
\[
\sqrt{-1}\left(\sum_{i=1}^r\frac{dz_i\wedge d \bar z_i}{|z_i|^2(\log\frac{1}{|z_i|})^2}
+\sum_{i=r+1}^{r+s} dz_i\wedge d\bar z_i\right).
\]
Since $\mathcal M$ is quasi-projective, it can be covered by finitely
many open sets of the form $(\Delta^*)^r\times\Delta^s$ ($r$ is
allowed to be zero). By Jost-Yau, we can glue the K\"ahler metrics of
the above form and get a global K\"ahler metric $\omega_P$ on
$\mathcal M$. The metric satisfies the following properties:

\begin{enumerate}
\item $\omega_P$  is complete;
\item The Ricci curvature of $\omega_P$ is bounded from below;
\item the volume of the metric $\omega_P$ is finite.
\end{enumerate}

Unlike the Weil-Petersson metric or Hodge metric, $\omega_P$ is not
intrinsically defined.

If we let $f$ in Theorem~\ref{yau} be the identity map from $\mathcal
M$ to itself, then using the Schwarz lemma, we get the following

\begin{lemma}\label{lem10}
 Let $\omega_H,\omega_P$ be the two metrics on $\mathcal M$. Then
 there is a constant $C$ such that
\[
\omega_H\leq C\omega_P.
\]
\end{lemma}

\qed

We remark that by ~\cite[Theorem A.1]{fl1}, even for the generalized
Hodge metric, the inequality in Lemma~\ref{lem10} is still valid. In
particular, this implies that the Hodge volumes and the Weil-Petersson
volume are all finite on a Calabi-Yau moduli space.

\section{Proof of the results.}

{\bf Proof of Theorem~\ref{main1}.}  Let $c_r(\omega_{WP})$ be the
$r$-th elementary polynomial of the curvature matrix of the
Weil-Petersson metric. We claim that
\begin{equation}\label{claim}
|c_r(\omega_{WP})|
\leq C\omega_{H}^r.
\end{equation}

The above inequality means that for any $v_1,\cdots, v_r\in T\mathcal
M$, we have
\[
|c_r(\omega_{WP})(v_1,\cdots,v_r,\bar v_1,\cdots,\bar v_r)|\leq C\prod_{i=1}^r||v_i||^2
\]
for some constant $C>0$, where the norm of the right hand side is with
respect to the metric $\omega_{H}$.

To prove the claim, first we choose a
normal coordinate system at $x\in \mathcal M$ such that  
\begin{equation}\label{dg}
g_{i\bar{j}}(x)=\delta_{ij},\quad dg_{i\bar{j}}(x)= 0.
\end{equation}

Let
\[ R_i^j = \sum_{kl}R_{ik\bar{l}}^jdz^k\wedge d\bar{z}^l, \quad \textrm{where 
$R^j_{ik\bar{l}}=g^{j\bar{p}}R_{i\bar{p}k\bar{l}}.$} \]
Then the $r$-th Chern class is given by
\begin{equation}\label{chern}
 c_r(\omega_{WP})=\frac{(-1)^r}{r!}\sum_{\tau\in S_r} sgn(\tau) R_{i_1}^{i_{\tau(1)}}\wedge \cdots \wedge R_{i_r}^{i_{\tau(r)}}, 
 \end{equation}
where $S_r$ is the symmetric group of the set $\{1,2,\ldots,r\}$.

We define

\[ h'_{\alpha\bar{\beta}} = \gd_{\alpha\beta}+\sum_\gamma(\nabla_\alpha \nabla_\gamma\Omega,\overline{\nabla}_\beta \overline{\nabla}_\gamma\Omega). \]
Then $(h_{\alpha\bar\beta}')$ defines a K\"ahler metric $\omega'$.
By~\cite[Proposition 2.8]{fl1} and Theorem~\ref{thm1},  we have
\[
\omega'\leq\omega_H.
\]
Thus in order to prove the claim, we only  need  to prove that
\[
|c_r(\omega_{WP})|\leq C(\omega')^r.
\]

Let $A_{ij}=\sum_k(\nabla_i \nabla _k\Omega,\overline{\nabla} _j\overline{\nabla} _k\Omega)$. 
 Since the matrix $(A_{ij})$ is Hermitian, after suitable
unitary change of basis, we can assume
\[ A_{ij}(p) = \left\{ \begin{array}{ll}
\lambda_i & \textrm{if $i=j$} \\
0 & \textrm{if $i\ne j$}.
\end{array}\right. \]		
Since $(A_{ij}(p))$ is positive-semidefinite, $\lambda_i \ge 0$, and we can write
\begin{equation}\label{dh}
h'_{i\bar{j}}(p) = \delta_{ij}(1+\lambda_i).
\end{equation}

Clearly, we have
\begin{equation}\label{deth}
(1+\lambda_{i_1}) \cdots (1+\lambda_{i_\ga}) \le \det h'
\end{equation}
for any $1\le i_1 < i_2<\ldots <i_\ga \le n$, where 
$\det h'=\det (h'_{\alpha\bar\beta})$. We assume that $v_i=\frac{\pa}{\pa t_{k_i}}$, then by~\eqref{chern}, we have
\[ |c_r(\omega_{WP})   (v_1,\cdots,v_r,\bar v_1,\cdots,\bar v_r)|\leq C\,{\rm Max}\,
|R_{j_1k_1\overline{\sigma(k_{1})}}^{i_1} \cdots R_{j_rk_r\overline{\sigma(k_r)}}^{i_r}| , \]
For fixed $i,j,k,l$, by the Cauchy-Schwartz inequality, we have 
\begin{eqnarray*}
\big|R^j_{ik\bar{l}}\big| &\le& |\gd_{ij}\gd_{kl}+\gd_{il}\gd_{kj}- (\nabla _i\nabla _k\Omega,\overline{\nabla} _j\overline{\nabla}_l\Omega)| \\
&\le& 2+ \sqrt{(\nabla _i\nabla _k\Omega,\overline{\nabla} _i\overline{\nabla} _k\Omega)(\nabla _j\nabla _l\Omega,\overline{\nabla} _j\overline{\nabla} _l\Omega)}  \\
&\le& 2+ \sqrt{h_{i\bar{i}}h_{j\bar{j}}}\leq 2  \sqrt{(1+\lambda_i)(1+\lambda_j)}.
\end{eqnarray*}
So we get
\[ \Big|R_{j_1k_1\bar{k}_{\sigma(1)}}^{i_1}\ldots R_{j_rk_r\bar{k}_{\sigma(r)}}^{i_r} \Big| \le
2^m\prod_{\ga=1}^m \Big(\sqrt{(1+\lambda_{k_i})(1+\lambda_{\sigma(k_i)}}) \; \Big). \]

The claim follows from the above inequality and ~\eqref{deth}.

Let $r_0,\cdots,r_t\geq 0$ such that $\sum r_i=m$. Then using~\eqref{claim},
\[
|c_{r_1}(\omega_{WP})\wedge\cdots \wedge c_{r_t}(\omega_{WP})\wedge \omega_{WP}^{r_0}|\leq C(\det h').
\]
By Lemma~\ref{lem10}, the left hand side of the above equation is integrable.
Theorem~\ref{main1} follows from the above inequality.

\qed

To prove Theorem~\ref{main2},  we first observe the following:

Let $E\rightarrow X$ be a holomorohic vector bundle over a compact manifold $X$. Let $h_0,h_1$ be two Hermitian metrics on the bundle. Let $R_0,R_1$ be the curvature tensors and let $\theta_0, \theta_1$ be the connection matrices. Let $f$ be an invariant polynomial. Then we have
\[
\int_X f(R_0,\cdots,R_0)=\int_X f(R_1,\cdots,R_1).
\]
In fact, 
\[
\int_X f(R_1,\cdots,R_1)-\int_X f(R_0,\cdots,R_0)
=\int_X\sum_{i=1}^k f(R_1,\cdots,R_1, R_1-R_0,R_0,\cdots,R_0).
\]
Since $R_1-R_0=\bar\pa(\theta_1-\theta_0)$,  we get
\begin{align*}
&
\int_X f(R_1,\cdots,R_1)-\int_X f(R_0,\cdots,R_0)\\
&
=\sum_{i=1}^k\int_X\bar\pa f(R_1,\cdots,R_1,\theta_1-\theta,R_0,\cdots,R_0)=0.
\end{align*}

A similar method can be used in the non-compact cases. The only
difference is that we need various estimates of the curvatures and the
connections near the infinity.

{\bf Proof of Theorem~\ref{main2}.}  As discussed in Section 3, up to
 a finite cover, we can assume that both $\mathcal M$ and
 $\overline{\mathcal M}$ are manifolds and $D=\overline{\mathcal
 M}\backslash{\mathcal M}$ is a divisor of normal crossings.
 Furthermore, by~\cite[Lemma 4.1]{ls-2}, we assume that the monodromy
 operators of the divisor $D$ are all unipotent. Since $\dim\mathcal
 M=2$, if $D_0$ denotes the smooth part of $D$, then the singular part
 $D\backslash D_0$ are the set of finite points. Let $F=F^k$ be a
 Hodge bundle and let $\bar F$ be the Schmid extension of the bundle
 across $D$.

  Let 
\[
D\backslash D_0=\{x_1,\cdots,x_s\}.
\]
Let $h$ be the Hermitian metric of $F$.
Let $U$ be a neighborhood of $D$ such that $x_i\notin\overline U$ for any $i$. Let $e_1,\cdots,e_k$ be a local holomorphic frame. Let $\langle e_i,e_j\rangle$ be the inner product induced from the metric $h$. Let $(z_1,z_2)$ be the local coordinates of $U$ such that $D\cap U=\{z_1=0\}$. Then by the Nilpotent orbit theorem of Schmid, we can write
\[
e_i=\exp(N\log 1/z_1) A_i(z_1,z_2), \quad 1\leq i\leq {\rm rank}\, F.
\]
It follows that the determinant of the metric matrix $\langle e_i,e_j\rangle$ can be expanded as 
\[
\det \langle e_i,e_j\rangle=a(z_1,z_2)(\log\frac{1}{r_1})^\alpha+\cdots,
\]
where $\cdots$ are the lower order terms and $a(z_1,z_2)$ is a real
analytic function of $(z_1,z_2)$. The zero set of $a(z_1,z_2)$ is
finite, and is independent to the choices of local frames. Let
$x_{s+1},\cdots,x_t$ be such kind of zeros on $D$. Let
$U_1,\cdots,U_t$ be neighborhoods of $x_i (1\leq i\leq t)$. Assume
that $U_i\cap U_j=\emptyset$ for $i\neq j$. These open sets are called
the neighborhoods of the first kind. Let $\{U_1,\cdots,U_{t+t_1}\}$ be
a cover of $D$.  The neighborhood $U_i (t<i\leq t+t_1)$ are called
neighborhoods of the second kind. Let $U_i$ be a neighborhood of the
second kind and let $(z_1,z_2)$ be the holomorphic coordinates. Let
$\theta, R$ be the connection and the curvature matrices of the metric
$h$. Let
\[
\theta=\theta_i dz_i,\qquad R=R_{i\bar j}dz_i\wedge d\bar z_j.
\]
Then there is a constant $C>0$ such that
\begin{align}\label{b1}
\begin{split}
& |\theta_1|\leq\frac{C}{r_1\log\frac{1}{r_1}}, \quad |\theta_2|\leq C;\\
&|R_{1\bar 1}|\leq\frac{C}{(r_1\log\frac{1}{r_1})^2}, \quad |R_{1\bar 2}|, |R_{2\bar 1}|\leq\frac{C}{r_1\log\frac{1}{r_1}},\quad  |R_{2\bar 2}|\leq C.
\end{split}
\end{align}

Let $h'$ be a Hermitian metric on $\overline F$ and let $\theta', R'$
be the corresponding connection and curvature matrices. We assume that
on the neighborhoods of the first kind, the metric is flat. That is,
$\theta'$ and $R'$ are identically zero on $U_i (1\leq i\leq t)$. Let
$U_1,\cdots, U_{t+t_1+t_2}$ be a cover of $\overline{\mathcal M}$ such
that
\begin{enumerate}
\item $U_i (1\leq i\leq t)$ are neighborhoods of the first kind;
\item $U_i (t\leq i\leq t+t_1)$ are neighborhoods of the second kind;
\item $D\cap (\bigcup_{i=t+t_1+1}^{t+t_1+t_2} \overline{U}_i)=\emptyset$.
\end{enumerate}

By~\cite[Theorem 3.1]{ls-2}, we know that for  any $\eps>0$, there is a cut-off function $\rho=\rho_\eps$ such that

\begin{enumerate}
\item
$0\leq\rho_\eps\leq 1$;
\item
For any open neighborhood $V$ of $D$ in $\overline{\mathcal  M}$, there is
$\eps>0$ such that ${\rm supp}(1-\rho_\eps)\subset V$;
\item For each $\eps>0$, there is a neighborhood
$V_1$ of $D$ such that $\rho_\eps|_{V_1}\equiv 0$;
\item $\rho_{\eps'}\geq\rho_\eps$ for $\eps'\leq \eps$;
\item There is a constant $C$, independent of $\eps$ such
that
\[
-C\omega_P\leq\sqrt{-1}\pa\bar\pa\rho_\eps\leq C\omega_P,\quad \left |\frac{\pa \rho}{\pa z_1}\right|\leq\frac{C}{r_1\log\frac{1}{r_1}}, \quad \left|\frac{\pa\rho}{\pa z_2}\right|\leq C.
\]
\end{enumerate}

Since
\[
R-R_0=\bar\pa(\theta-\theta_0),
\]
we have
\begin{equation}\label{s1}
\int_{\mathcal M} \rho(f(R, R)-f(R_0, R_0))=-\int_{\mathcal M}\bar\pa\rho\wedge f(R,\theta-\theta_0)-\int_{\mathcal M}\bar\pa\rho\wedge f(R_0,\theta-
\theta_0),
\end{equation}
where $f$ is an invariant quadratic polynomial.
For $\eps>0$ small enough, 
On $U_i (t+t_1<i\leq t+t_1+t_2)$, $\rho\equiv 1$. So we have
\begin{align*}&
\left |\int_{\mathcal M} \rho(f(R, R)-f(R_0, R_0))\right|\\&\leq\sum_{i=1}^{t+t_1}\left(\left|\int_{U_i}\bar\pa\rho\wedge f(R,\theta-\theta_0)\right|+\left|\int_{U_i}\bar\pa\rho\wedge f(R_0,\theta-
\theta_0)
\right|.\right)
\end{align*}
We shall prove that the right hand side of the above goes to zero as
$\eps\rightarrow 0$. If $U_i$ is a neighborhood of the second kind,
then since $\theta_0$ is bounded, by ~\eqref{b1} and the definition of
$\rho$, we know that
\[
|\bar\pa\rho\wedge f(R,\theta-\theta_0)|+|\bar\pa\rho\wedge f(R_0,\theta-
\theta_0)|\leq \frac{C}{r_1^2(\log\frac{1}{r_1})^2}|dz_1\wedge d\bar z_1\wedge dz_2\wedge d\bar z_2|,
\]
which is integrable. Thus we have
\begin{align*}
&
\lim_{\eps\rightarrow 0}\left|\int_{U_i}\bar\pa\rho\wedge f(R,\theta-\theta_0)\right|+\left|\int_{U_i}\bar\pa\rho\wedge f(R_0,\theta-
\theta_0)
\right|\\
&\leq C\int_{{\rm supp}\,\nabla \rho}({r_1^2(\log\frac{1}{r_1})^2})^{-1}|dz_1\wedge d\bar z_1\wedge dz_2\wedge d\bar z_2|=0
\end{align*}
for $t<i\leq t+t_1$. 
If $U_i$ is a neighborhood of the first kind, then we have
\[
\int_{U_i}\bar\pa\rho\wedge f(R_0,\theta-
\theta_0)
=0
\]
because $R_0\equiv 0$ on $U_i$. 
On the other hand, we have
\[
\int_{U_i}\bar\pa\rho\wedge f(R,\theta-\theta_0)=\int_{U_i}\bar\pa\rho\wedge f(R,\theta).
\]
From~\cite[Proposition 5.22]{cks1}, we know that there is a gauge transform $e$ such that we have
\begin{align*}
& |Ad(e)\theta_1|\leq\frac{C}{r_1\log\frac{1}{r_1}}, \quad |Ad(e)\theta_2|\leq C;\\
&|Ad(e)R_{1\bar 1}|\leq\frac{C}{(r_1\log\frac{1}{r_1})^2}, \quad |Ad(e)R_{1\bar 2}|, |Ad(e)R_{2\bar 1}|\leq\frac{C}{r_1\log\frac{1}{r_1}},\quad  |Ad(e)R_{2\bar 2}|\leq C.
\end{align*}
Since $f$ is an invariant polynomial, using the transform, we have
\begin{align*}
&
\left|\int_{U_i}\bar\pa\rho\wedge f(R,\theta)\right|=\left|
\int_{U_i}\bar\pa\rho\wedge f(Ad(e)R,Ad(e)\theta)\right|\\
\leq &
C\int_{{\rm supp}\,\nabla \rho}({r_1^2(\log\frac{1}{r_1})^2})^{-1}|dz_1\wedge d\bar z_1\wedge dz_2\wedge d\bar z_2|\rightarrow 0.
\end{align*}

Thus from~\eqref{s1}, we proved that
\[
\lim_{\eps\rightarrow 0} \int_{\mathcal M} \rho(f(R, R)-f(R_0, R_0))=0.
\]

By the Gauss-Bonnet-Chern theorem, 
\[
\int_{\mathcal M} f(R_0, R_0)
\]
is an integer. Thus
\[
\int_{\mathcal M} f(R, R)
\]
is also an integer.

\qed

\begin{remark} Theorem~\ref{main2} is 
a rationality result of the Hodge bundles. However, by
Proposition~\ref{pp1}, the Weil-Petersson metric is the quotient of
the Hermitian metrics on the Hodge bundles $H^{n-1,1}$ and $H^{n,0}$,
respectively. Thus the corresponding result for the Weil-Petersson
metric is also valid.
\end{remark}

\bibliographystyle{abbrv} 
\bibliography{new051007,unp051007,local}

\end{document}